\documentclass[12pt]{article} 

\usepackage{latexsym, 
amscd, 
graphicx, epsfig, amssymb}

\newtheorem{thm}{Theorem}[section]
\newtheorem{defn}[thm]{Definition}

\newtheorem{rema}[thm]{Remark}

\newcommand{\halmos}{\rule{1ex}{1.4ex}}

\newcommand{\nn}{\nonumber \\}

 \newcommand{\res}{\mbox{\rm Res}}
\renewcommand{\hom}{\mbox{\rm Hom}}
 \newcommand{\pf}{{\it Proof.}\hspace{2ex}}
 
 \newcommand{\epfv}{\hspace*{\fill}\mbox{$\halmos$}\vspace{1em}}

\newcommand{\wt}{\mbox{\rm wt}\ }

\newcommand{\lbar}{\bigg\vert}

\newcommand{\injto}{\hookrightarrow}  
\newcommand{\one}{\mathbf{1}}
\renewcommand{\th}{\theta}

\title{ {\bf Braided tensor categories and extensions 
of vertex operator algebras} }
\date{}
\author{Yi-Zhi Huang, Alexander Kirillov Jr. and James Lepowsky}

\begin{document}

\bibliographystyle{alpha}
\maketitle

\begin{abstract} 
Let $V$ be a vertex operator algebra satisfying suitable conditions
such that in particular its module category has a natural vertex
tensor category structure, and consequently, a natural braided tensor
category structure.  We prove that the notions of extension (i.e.,
enlargement) of $V$ and of commutative associative algebra, with
uniqueness of unit and with trivial twist, in the braided tensor
category of $V$-modules are equivalent.
\end{abstract}

\renewcommand{\theequation}{\thesection.\arabic{equation}}
\renewcommand{\thethm}{\thesection.\arabic{thm}}
\setcounter{equation}{0}
\setcounter{thm}{0}

\section{Introduction}

In the constructions and applications of vertex operator algebras,
such algebras obtained by adjoining suitable modules to a smaller
vertex operator algebra play a significant role, as in \cite{FLM}, for
example.

Given a vertex operator algebra $V$ satisfying suitable conditions, a
natural vertex tensor category structure, and as a consequence, a
natural braided tensor category structure, on the category of
$V$-modules was constructed by the first and the third authors in
\cite{HL1}--\cite{HL7}, \cite{H1} and \cite{H4}.  Such braided tensor
category structure on the category of $V$-modules (including even the
construction of suitable tensor product modules) had been conjectured
to exist in physics and mathematics starting from the work of Moore
and Seiberg \cite{MS}.  Assuming the existence of such braided tensor
category structure, it was proposed in \cite{Ho}, \cite{FFFS} and
\cite{KO} that certain extensions of a vertex operator algebra $V$ (by
which we mean vertex operator algebras containing $V$ as a subalgebra,
in this paper) should be related to commutative associative algebras
in the braided tensor category of $V$-modules. In \cite{KO}, Ostrik
and the second author formulated a relation between the extensions of
$V$ and the commutative associative algebras in the category of
$V$-modules.  Assuming this relation, in \cite{KO} they classified the
extensions of the vertex operator algebras associated to standard
(integrable highest weight) modules for the affine Lie algebra
$\widehat{\mathfrak{sl}}_2$, relating these extensions to the
representation theory of $U_q({\mathfrak{sl}}_2)$ and obtaining a
$q$-analogue of the McKay correspondence.

In the present paper, we establish close variants of this relation
formulated in \cite{KO} between the extensions of $V$ and the
commutative associative algebras in the category of $V$-modules.  This
requires the theory of vertex tensor categories.  There are many
classes of vertex operator algebras whose suitable module categories
are vertex tensor categories (see the Appendix below), and for these
module categories, the conclusions of the theorems in the present
paper hold, under suitable additional conditions.  These classes of
vertex operator algebras include in particular the vertex operator
algebras associated with positive definite even lattices; the vertex
operator algebras associated with affine Lie algebras and positive
integral levels; the ``minimal series'' of vertex operator algebras
associated with the Virasoro algebra; Frenkel, Lepowsky and Meurman's
moonshine module vertex operator algebra $V^{\natural}$; the fixed
point vertex operator subalgebra of $V^{\natural}$ under the standard
involution; and the ``minimal series'' of vertex operator
superalgebras (suitably generalized vertex operator algebras)
associated with the Neveu-Schwarz superalgebra and also the ``unitary
series'' of vertex operator superalgebras associated with the $N=2$
superconformal algebra.  Since the construction of the necessary
vertex tensor category structure is substantial, the proofs of these
conclusions for these module categories are highly nontrivial.

More specifically, in contrast with quantum group theory, in which the
associativity isomorphisms are trivial, the construction of the
natural associativity isomorphisms for the vertex tensor category
structures and the braided tensor category structures, and the proofs
of their properties, form the hardest part of the work
\cite{HL1}--\cite{HL7}, \cite{H1}, \cite{H4}.  (See Section 2 below
for a sketch of the necessary structures.)  In particular, using
strict tensor categories equivalent to the original ones would only
make the necessary constructions and proofs more difficult, and this
is true even in the monodromy-free case when the intertwining maps
(induced by intertwining operators) are single-valued.

In our theorems, we always assume that the module category for the
underlying vertex operator algebra has a natural vertex tensor
category structure rather than just a braided tensor category
structure, since our proofs are based on this structure.  In our first
theorem, we do not require that the braided tensor category be rigid.
We obtain different variants of the relation formulated in \cite{KO}.

The first version of the present paper was written in 2002, and
mathematically, the present paper is essentially the same as the 2002
version.  The results have been generalized by Kong and the first
author in \cite{HK} to a relation between open-string vertex algebras
containing $V$ as a subalgebra and associative algebras in the
category of $V$-modules, and by Kong in \cite{K} to a relation between
conformal full field algebras over $V^{L}\otimes V^{R}$, for suitable
vertex operator algebras $V^{L}$ and $V^{R}$ equipped with
nondegenerate invariant bilinear forms, and commutative Frobenius
algebras with trivial twists in the category of $V^L\otimes
V^R$-modules.

In the next section, we review the relevant aspects of the vertex
tensor category structure (\cite{HL1}--\cite{HL7}, \cite{H1},
\cite{H4}).  In Section 3 we formulate and prove our theorems,
including the result that the notion of vertex operator algebra
extension of a vertex operator algebra satisfying certain conditions
is equivalent to the notion of commutative associative algebra, with
uniqueness of unit and with trivial twist, in the braided tensor
category of $V$-modules.  
In the Appendix, for the reader's convenience we recall general
theorems stating that under suitable conditions, the module category
for a vertex operator algebra has a natural vertex tensor category
structure and in particular braided tensor category structure.

\paragraph{Acknowledgment} 
The first and third author were partially supported 
by NSF grants DMS-0070800 and DMS-0401302. 

\renewcommand{\theequation}{\thesection.\arabic{equation}}
\renewcommand{\thethm}{\thesection.\arabic{thm}}
\setcounter{equation}{0}
\setcounter{thm}{0}

\section{Review of the vertex tensor category theory}

Here we briefly review what we will need concerning vertex tensor
category structure for the proofs of the main theorems of the present
paper.  This entails the vertex tensor category structure itself,
including the resulting braided tensor category structure, on the
category of modules for a vertex operator algebra $V$ satisfying
suitable conditions; this structure was constructed in
\cite{HL1}--\cite{HL7}, \cite{H1}, \cite{H4}, whose notation and
terminology we shall be using.  For suitable sufficient conditions,
see the Appendix.

We take as our precise notions of vertex operator algebra (over
$\mathbb{C}$) and of module (over $\mathbb{C}$) for a vertex operator
algebra the definitions given in \cite{HL3}. In particular, the
central charge (or rank) of a vertex operator algebra is allowed to be
in $\mathbb{C}$ and the weight-grading of a module is by $\mathbb{C}$
(that is, the (conformal) weights of homogeneous elements are complex
numbers).

We shall use the notation $Y$ for the vertex operator map of the
vertex operator algebra $V$, $\mathbf{1}$ for the vacuum vector,
$\omega$ for the conformal element (generating the Virasoro algebra),
$L(n)$ ($n \in \mathbb{Z}$) for the Virasoro algebra operators, and
$V_{(n)}$ ($n \in \mathbb{Z}$) for the $L(0)$-eigenspaces (weight
spaces) of $V$, as well as the analogous relevant notation for
$V$-modules.  By definition, a vertex operator subalgebra of a vertex
operator algebra has the same conformal element.  A simple vertex
operator algebra is one that is irreducible as a module for itself.

We shall need the $P(z)$-tensor product $W_{1}\boxtimes_{P(z)}W_{2}$
of $V$-modules $W_{1}$ and $W_{2}$, for $z\in \mathbb{C} \setminus
\{0\}$.  There is a tensor product bifunctor associated to any sphere
with three tubes, but in the present paper the $P(z)$-tensor product
is enough.  The $V$-module $W_{1}\boxtimes_{P(z)}W_{2}$ is not based
on the tensor product vector space $W_{1}\otimes W_{2}$.

One important feature of the $P(z)$-tensor product
$W_{1}\boxtimes_{P(z)}W_{2}$ is a natural isomorphism between the
space of module maps from $W_{1}\boxtimes_{P(z)} W_{2}$ to a third
module $W_{3}$ and the space of intertwining operators of type
${W_{3}\choose W_{1} W_{2}}$ (see \cite{HL3}, \cite{HL5}).  (This
isomorphism is based on the choice of the branch of $\log z$ such that
$0\le \arg z<2\pi$; different choices give equivalent theories.)  In
particular, there exists a canonical intertwining operator
$\mathcal{Y}(\cdot,x)$ ($x$ a formal variable) of type
${W_{1}\boxtimes_{P(z)}W_{2} \choose W_{1}\;W_{2}}$ corresponding to
the identity module map from $W_{1}\boxtimes_{P(z)} W_{2}$ to itself.
For $w_{1}\in W_{1}$ and $w_{2}\in W_{2}$, one has a $P(z)$-tensor
product element
\[
w_{1}\boxtimes_{P(z)} w_{2} =\mathcal{Y}(w_{1}, z)w_{2} \in
\overline{W_{1}\boxtimes_{P(z)}W_{2}}
\]
of $w_{1}$ and $w_{2}$, with $z$ replacing $x$ according to the branch
choice, giving a $P(z)$-intertwining map $\mathcal{Y}(\cdot,z)$ (as
opposed to intertwining operator), where
$\overline{W_{1}\boxtimes_{P(z)}W_{2}}$ is the formal completion of
$W_{1}\boxtimes_{P(z)}W_{2}$.  Moreover, the homogeneous components of
$w_{1}\boxtimes_{P(z)} w_{2}$ for all $w_{1}\in W_{1}$ and $w_{2}\in
W_{2}$ span the module $W_{1}\boxtimes_{P(z)}W_{2}$.

The existence of these tensor product elements is another important
feature of the tensor category theory.  Such elements provide a
powerful method for proving the necessary theorems in the tensor
category theory (see in particular \cite{H1}): One first proves
certain results about these tensor product elements.  Since the
homogeneous components of these elements span the tensor product
modules, one obtains the desired results.  (Worth mentioning here is
the subtlety that the set of all tensor product elements has almost
no intersection with the $P(z)$-tensor product module; generically,
only very special elements like $\mathbf{1}\boxtimes_{P(z)} w$ lie in
$V\boxtimes_{P(z)}W$, $W$ a $V$-module.)

To prove the necessary results involving associativity and coherence,
one also needs tensor product elements of more than two elements.
Here we briefly describe the tensor product element of three elements
(see \cite{H1} for details).  Let $W_{1}, W_{2}, W_{3}$ be $V$-modules
and let $w_{1}\in W_{1}$, $w_{2}\in W_{2}$ and $w_{3}\in W_{3}$.  Let
$z_{1}$ and $z_{2}$ be nonzero complex numbers.  Since in general
$w_{2}\boxtimes_{P(z_{2})} w_{3}$ does not lie in
$W_{2}\boxtimes_{P(z_{2})}W_{3}$, one cannot define
$w_{1}\boxtimes_{P(z_{1})}(w_{2}\boxtimes_{P(z_{2})} w_{3})$ simply as
the $P(z_{1})$-tensor product element of $w_{1}$ and
$w_{2}\boxtimes_{P(z_{2})} w_{3}$.  Such triple tensor product
elements do not exist in general.  But using the assumed convergence
property for the product of intertwining operators, the series
$\sum_{n\in \mathbb{R}}
w_{1}\boxtimes_{P(z_{1})}\pi_{n}(w_{2}\boxtimes_{P(z_{2})} w_{3})$
($\pi_{n}$ being the projection map from a graded space to the
subspace of weight $n$) is weakly absolutely convergent when
$|z_{1}|>|z_{2}|>0$ and the sum, which one writes as
$w_{1}\boxtimes_{P(z_{1})}(w_{2}\boxtimes_{P(z_{2})} w_{3})$, lies in
$\overline{W_{1}\boxtimes_{P(z_{1})}(W_{2} \boxtimes_{P(z_{2})}
W_{3})}$.  Similarly, when $|z_{2}|>|z_{1}-z_{2}|>0$, one has
$(w_{1}\boxtimes_{P(z_{1}-z_{2})}w_{2})\boxtimes_{P(z_{2})} w_{3}\in
\overline{(W_{1}\boxtimes_{P(z_{1}-z_{2})}W_{2})\boxtimes_{P(z_{2})}
W_{3}}$. The homogeneous components of these triple tensor product
elements also span the corresponding triple tensor product modules.
More generally, one has
\[
\mathcal{Y}_{1}(w_{1}, z_{1})\mathcal{Y}_{2}(w_{2}, z_{2})w_{3}
\]
and 
\[
\mathcal{Y}^{1}(\mathcal{Y}^{2}(w_{1}, z_{1}-z_{2})w_{2}, z_{2})w_{3}
\]
in the formal completions of suitable $V$-modules for suitable
intertwining operators $\mathcal{Y}_{1}$, $\mathcal{Y}_{2}$,
$\mathcal{Y}^{1}$ and $\mathcal{Y}^{2}$, elements $w_{1}$, $w_{2}$ and
$w_{3}$ of suitable $V$-modules, and complex numbers $z_{1}$ and
$z_{2}$ satisfying $|z_{1}|>|z_{2}|>0$ and $|z_{2}|>|z_{1}-z_{2}|>0$,
respectively.

The $P(z)$-tensor product functors and elements are needed for the
vertex tensor category structure and for the proof of the main theorems
of the present paper.  For the braided tensor category structure, one
needs to specify one particular tensor product bifunctor.  One can in
fact take this tensor product for the braided tensor category to be
any particular $P(z)$-tensor product, but one can and does choose the
$P(1)$-tensor product $\boxtimes_{P(1)}$ (where $z=1$), and one also
denotes $\boxtimes_{P(1)}$ by $\boxtimes$.

The unit object of the vertex tensor category and of the braided
tensor category is the vertex operator algebra $V$.  The left unit
isomorphism
\[
l_{W}: V\boxtimes W \to W
\]
is characterized
by
\[
l_{W}(\mathbf{1}\boxtimes w)=w
\]
for $w\in W$.  The right unit isomorphism
\[
r_{W}: W\boxtimes V \to W
\]
is characterized by
\[
\overline{r_{W}}(w\boxtimes
\mathbf{1})=e^{L(-1)} w
\]
for $w\in W$.

One needs natural parallel transport isomorphisms.  Let $W_{1}$ and
$W_{2}$ be $V$-modules.  Given a path $\gamma$ in $\mathbb{C}\setminus
\{0\}$ from $z_{1}\in \mathbb{C}\setminus \{0\}$ to $z_{2}\in
\mathbb{C}\setminus \{0\}$, the parallel transport isomorphism
$\mathcal{T}_{\gamma}: W_{1}\boxtimes_{P(z_{1})}W_{2}\to
W_{1}\boxtimes_{P(z_{2})}W_{2}$ is defined as follows: Let
$\mathcal{Y}$ be the intertwining operator corresponding to the
intertwining map $\boxtimes_{P(z_{2})}$ and let $l(z_{1})$ be the
value of the logarithm of $z_{1}$ determined uniquely by $\log z_{2}$
(satisfying $0\le \arg z_{2} <2\pi$) and the path $\gamma$. Then
$\mathcal{T}_{\gamma}$ is characterized by
\[
\overline{\mathcal{T}}_{\gamma}(w_{1}
\boxtimes_{P(z_{1})}w_{2})=\mathcal{Y}(w_{1}, e^{l(z_{1})})w_{2}
\]
for $w_{1}\in W_{1}$ and $w_{2}\in W_{2}$, where
$\overline{\mathcal{T}}_{\gamma}$ is the natural extension of
$\mathcal{T}_{\gamma}$ to the formal completion
$\overline{W_{1}\boxtimes_{P(z_{1})} W_{2}}$ of
$W_{1}\boxtimes_{P(z_{1})} W_{2}$.  The parallel transport isomorphism
depends only on the homotopy class of $\gamma$.

The braiding isomorphism (also called the commutativity isomorphism)
\[
\mathcal{R}_{W_{1}W_{2}}: W_{1}\boxtimes W_{2}\to W_{2}\boxtimes W_{1}
\]
is characterized as follows: Let $\gamma_{1}^{-}$ be a path from $-1$
to $1$ in the closed upper half plane with $0$ deleted, with
$\mathcal{T}_{\gamma_{1}^{-}}$ the corresponding parallel transport
isomorphism. Then
\begin{equation}
\overline{\mathcal{R}}_{W_{1}W_{2}}(w_{1}\boxtimes w_{2})=e^{L(-1)}
\overline{\mathcal{T}}_{\gamma_{1}^{-}}
(w_{2}\boxtimes_{P(-1)} w_{1})\label{commu-iso-1}
\end{equation}
for $w_{1}\in W_{1}$, $w_{2}\in W_{2}$. One can prove that 
\[
(\mathcal{R}_{W_{2}W_{1}}\circ \mathcal{R}_{W_{1}W_{2}})(w_{1}\boxtimes w_{2})
=\mathcal{Y}(w_{1}, x)w_{2}|_{x^{n}=e^{2\pi n i},\; n\in\mathbb{R}},
\]
where $w_{1}\in W_{1}$, $w_{2}\in W_{2}$ and $\mathcal{Y}$ is the
intertwining operator corresponding to the identity map from
$W_{1}\boxtimes W_{2}$ to itself.

Given complex numbers $z_{1}$ and $z_{2}$ satisfying
$|z_{1}|>|z_{2}|>|z_{1}-z_{2}|>0$ and $V$-modules $W_{1}$, $W_{2}$ and
$W_{3}$, one is able to construct an associativity isomorphism
\[
\mathcal{A}^{P(z_{1}-z_{2}), P(z_{2})}_{P(z_{1}), P(z_{2})}: 
W_{1}\boxtimes_{P(z_{1})} (W_{2}\boxtimes_{P(z_{2})} W_{3})\to 
(W_{1}\boxtimes_{P(z_{1}-z_{2})} W_{2})\boxtimes_{P(z_{2})} W_{3}
\]
characterized by 
\[
\overline{\mathcal{A}^{P(z_{1}-z_{2}), P(z_{2})}_{P(z_{1}), P(z_{2})}}
(w_{1}\boxtimes_{P(z_{1})} (w_{2}\boxtimes_{P(z_{2})} w_{3}))
=(w_{1}\boxtimes_{P(z_{1}-z_{2})} w_{2})\boxtimes_{P(z_{2})} 
w_{3}.
\]
To obtain the associativity isomorphism
\[
\mathcal{A}_{W_{1}, W_{2}, W_{3}}: 
W_{1}\boxtimes (W_{2}\boxtimes W_{3})\to 
(W_{1}\boxtimes W_{2})\boxtimes W_{3}
\]
for the braided tensor category structure, one needs some particular
parallel transport isomorphisms.  Let $z_{1}$ and $z_{2}$ be real
numbers satisfying $z_{1}>z_{2}>z_{1}-z_{2} > 0$.  Let $\gamma_{1}$
and $\gamma_{2}$ be paths in the real line with $0$ deleted from $1$
to $z_{1}$ and $z_{2}$, respectively, and $\gamma_{3}$ and
$\gamma_{4}$ be paths in the real line with $0$ deleted from $z_{2}$
and $z_{1}-z_{2}$ to $1$, respectively.  Then the associativity
isomorphism is given by
\begin{equation}
\mathcal{A}_{W_{1}, W_{2}, W_{3}}
=\mathcal{T}_{\gamma_{3}}\circ (\mathcal{T}_{\gamma_{4}}
\boxtimes_{P(z_{2})} 1_{W_3})\circ 
\mathcal{A}^{P(z_{1}-z_{2}), P(z_{2})}_{P(z_{1}), P(z_{2})}\circ
(1_{W_1} \boxtimes_{P(z_{1})} 
\mathcal{T}_{\gamma_{2}})\circ \mathcal{T}_{\gamma_{1}}.\label{mu-assoc5}
\end{equation}

The characterizations of the commutativity and associativity
isomorphisms above make the proof of the coherence properties
straightforward.  Here we sketch the proof of the commutativity of the
pentagon diagram. This commutativity states that the following
identity for maps from $W_{1}\boxtimes (W_{2}\boxtimes (W_{3}
\boxtimes W_{4}))$ to $((W_{1}\boxtimes W_{2})\boxtimes W_{3})
\boxtimes W_{4}$ holds:
\begin{eqnarray}\label{pentagon}
\lefteqn{\mathcal{A}_{W_{1}\boxtimes W_{2}, W_{3}, W_{4}}\circ
 \mathcal{A}_{W_{1}, W_{2}, W_{3}\boxtimes W_{4}}}\nn
&&=(\mathcal{A}_{W_{1}, W_{2}, W_{3}}\boxtimes 1_{W_{4}}) \circ 
\mathcal{A}_{W_{1}, W_{2}\boxtimes W_{3}, W_{4}}
\circ (1_{W_{1}}\boxtimes \mathcal{A}_{W_{2}, W_{3}, W_{4}}).
\end{eqnarray}
By analogy with the construction of the tensor product elements
$w_{1}\boxtimes_{P(z_{1})} (w_{2}\boxtimes_{P(z_{2})} w_{3})$, one can
also construct tensor product elements of four elements, including
\begin{equation}\label{element1}
w_{1}\boxtimes_{P(z_{1})} (w_{2}\boxtimes_{P(z_{2})}
(w_{3} \boxtimes_{P(z_{3})} w_{4}))
\end{equation}
and
\begin{equation}\label{element2}
((w_{1}\boxtimes_{P(z_{1}-z_{2})} w_{2})\boxtimes_{P(z_{2}-z_{3})}
w_{3}) \boxtimes_{P(z_{3})} w_{4}))
\end{equation}
for $w_{1}\in W_{1}$, $w_{2}\in W_{2}$, $w_{3}\in W_{3}$, $w_{4}\in W_{4}$,
when 
\[
|z_{1}|>|z_{2}|>|z_{3}|>|z_{1}-z_{3}|>|z_{2}-z_{3}|>|z_{1}-z_{2}|>0.
\]
Since $\overline{\mathcal{A}^{P(z_{1}-z_{2}), P(z_{2})}_{P(z_{1}),
P(z_{2})}}$ maps $w_{1}\boxtimes_{P(z_{1})} (w_{2}\boxtimes_{P(z_{2})}
w_{3})$ to $(w_{1}\boxtimes_{P(z_{1}-z_{2})}
w_{2})\boxtimes_{P(z_{2})} w_{3}$, the formal extensions of both
\[
\mathcal{A}^{P(z_{2}-z_{3}), P(z_{3})}_{P(z_{2}), P(z_{3})}\circ
\mathcal{A}^{P(z_{1}-z_{2}), P(z_{2})}_{P(z_{1}), P(z_{2})}
\]
and 
\[
(\mathcal{A}^{P(z_{1}-z_{2}), 
P(z_{2}-z_{3})}_{P(z_{1}-z_{3}), P(z_{2}-z_{3})}
\boxtimes_{P(z_{3})} 1_{W_{4}})\circ
\mathcal{A}^{P(z_{1}-z_{3}), P(z_{3})}_{P(z_{1}), P(z_{3})}
\circ (1_{W_{1}}\boxtimes_{P(z_{1})} 
\mathcal{A}^{P(z_{2}-z_{3}), P(z_{3})}_{P(z_{2}), P(z_{3})})
\]
map (\ref{element1}) to (\ref{element2}).
Thus these two maps are the same:
\begin{eqnarray}\label{pentagon2}
\lefteqn{\mathcal{A}^{P(z_{2}-z_{3}), P(z_{3})}_{P(z_{2}), P(z_{3})}\circ
\mathcal{A}^{P(z_{1}-z_{2}), P(z_{2})}_{P(z_{1}), P(z_{2})}}\nn
&&=
(\mathcal{A}^{P(z_{1}-z_{2}), 
P(z_{2}-z_{3})}_{P(z_{1}-z_{3}), P(z_{2}-z_{3})}
\boxtimes_{P(z_{3})} 1_{W_{4}})\circ 
\mathcal{A}^{P(z_{1}-z_{3}), P(z_{3})}_{P(z_{1}), P(z_{3})}\nn
&&\quad\quad\quad
\circ (1_{W_{1}}\boxtimes_{P(z_{1})} 
\mathcal{A}^{P(z_{2}-z_{3}), P(z_{3})}_{P(z_{2}), P(z_{3})}).
\end{eqnarray}
Now one chooses $z_{1}$ and $z_{2}$ to be real numbers satisfying
$z_{1}>z_{2}>z_{1}-z_{2}>0$.  From the characterization above of the
associativity isomorphism for the braided tensor category structure
and the fact that the parallel transport isomorphisms involved are
associated to paths in the real line, the parallel transport
isomorphisms all cancel, and (\ref{pentagon}) follows from
(\ref{pentagon2}).

For any $V$-module $W$, one also has a twist $\theta_{W}\colon W\to
W$, given by $\theta_{W}(w)=e^{2\pi iL(0)}w$ for $w\in W$.  It is
immediate from the definition that if $W$ is simple, then
\[
\theta_W=e^{2\pi i \Delta_W},
\]
where $\Delta_W$ is the conformal dimension of $W$, that is, the
lowest eigenvalue of the weight operator $L(0)$ on $W$.  The twist
$\theta_W$ satisfies the following conditions:
\begin{eqnarray}
\theta_V &=& 1_V,\label{theta-1}\\
\theta_{W_1\boxtimes W_2}&=&\mathcal{R}_{W_2
W_1}\circ \mathcal{R}_{W_1 W_2}\circ 
(\theta_{W_1}\boxtimes\theta_{W_2}),\label{theta-2}\\
\theta_{W'}&=&(\theta_W)',\label{theta-3}
\end{eqnarray}
where $W'$ is the $V$-module contragredient to $W$, and for $f\colon
W\to W$ one denotes by $f'\colon W'\to W'$ the adjoint operator (see
\cite{FHL}).

The first identity (\ref{theta-1}) is immediate from the definitions.
For (\ref{theta-2}), note that for $w_{1}\in W_{1}$ and $w_{2}\in
W_{2}$, $w_{1}\boxtimes w_{2}=\mathcal{Y}(w_{1}, 1)w_{2}$ and
$\mathcal{R}_{W_2 W_1}\mathcal{R}_{W_1 W_2}(w_{1}\boxtimes
w_{2})=\mathcal{Y}(w_{1}, e^{2\pi i})w_{2}$ where $\mathcal{Y}$ is the
corresponding intertwining operator (see above). Thus
\[
\overline{\mathcal{R}_{W_2
W_1}\circ \mathcal{R}_{W_1 W_2}\circ 
(\theta_{W_1}\boxtimes\theta_{W_2})}(w_{1}\boxtimes w_{2})
=\mathcal{Y}(e^{2\pi i L(0)}w_1, e^{2\pi i })e^{2\pi i L(0)}w_2,
\]
and (\ref{theta-2}) follows from basic property
\[
e^{2\pi i L(0)}\mathcal{Y}(w_1, z)w_2=\mathcal{Y}(e^{2\pi iL(0)}w_1, 
e^{2\pi i})e^{2\pi iL(0)}w_2.
\]
The last identity (\ref{theta-3}) is immediate because $L(0)'=L(0)$.

If the tensor category of $V$-modules is rigid\footnote{After the
first version of the present paper was written, the rigidity of the
tensor category of $V$-modules was proved by the first author when $V$
satisfies certain stronger conditions.  See the Appendix for more
discussion and the reference.}, then $\theta$ defines a structure of a
balanced rigid category on this category.

\renewcommand{\theequation}{\thesection.\arabic{equation}}
\renewcommand{\thethm}{\thesection.\arabic{thm}}
\setcounter{equation}{0}
\setcounter{thm}{0}

\section{The main theorems}

\begin{defn}\label{dalg}
{\em Let $\mathcal{C}$ be a braided tensor category.  A {\em
commutative associative algebra} $A$ {\em in} $\mathcal{C}$ (or
$\mathcal{C}$-{\em algebra} for short) is an object $A$ of
$\mathcal{C}$ together with morphisms 
\[
\mu\colon A\otimes A\to A
\]
and
\[
\iota_A\colon \one_\mathcal{C} \injto A
\]
such that the following
conditions hold:
\begin{enumerate}

\item Associativity:
\[
\mu \circ(\mu \otimes 1_A)\circ \mathcal{A}=\mu\circ (1_A \otimes
\mu): A \otimes (A \otimes A) \longrightarrow A,
\]
where $\mathcal{A}$ is the associativity isomorphism
\[
\mathcal{A}:A\otimes (A\otimes A) \longrightarrow (A\otimes A)\otimes A.
\]

\item Commutativity:
\[
\mu \circ \mathcal{R}=\mu, 
\]
where $\mathcal{R}$ is the braiding isomorphism
\[
\mathcal{R}:A\otimes A \longrightarrow A\otimes A.
\]

\item Unit:
\[
\mu \circ (\iota_A \otimes 1_A)\circ l^{-1}_{A}=1_A,
\]
where
\[
l_{A}: \one_{\mathcal{C}}\otimes A\to A
\]
is the left unit isomorphism.
\end{enumerate}
Such an algebra $A$ is called {\em haploid} if it satisfies the
following additional condition:
\begin{enumerate}
\setcounter{enumi}{3}
\item Uniqueness of unit:
\[
\dim \hom_{\mathcal{C}}(\one_\mathcal{C}, A)=1.
\]

\end{enumerate}
}

\end{defn}

(Since we are using the notation $\mathbf{1}$ for the vacuum vector of
$V$, we use the notation $\one_\mathcal{C}$ for the unit object of
$\mathcal{C}$.)

Recall that a balancing in a rigid braided tensor category is a
natural isomorphism from the identity functor to the double-dual
functor satisfying standard conditions.  This is equivalent to
defining a natural transformation, called a twist, from the identity
functor on $\mathcal{C}$ to itself satisfying standard conditions.
Note that in the case of module categories for vertex operator
algebras, we have a twist $\theta_{W}$ on any module $W$ even if the
category is not rigid (see Section 2). Below we shall use this twist
$\theta_{W}$.

\begin{thm}\label{main}
Let $V$ be a vertex operator algebra such that its module category
$\mathcal{C}$ has a natural vertex tensor category structure described
in Section 2 (see the Appendix for the appropriate conditions on
$V$). We assume in addition that $V$ is simple; every $V$-module is
completely reducible; for any irreducible $V$-module $W$ the weights
of elements of $W$ are nonnegative real numbers and $W_{(0)}\ne 0$ if and only if
$W$ is equivalent to $V$ viewed as a $V$-module; and $\dim
V_{(0)}=1$. Then the following two notions are equivalent:

\begin{enumerate}
\item An extension $V_{e}\supset V$ such that $\dim (V_{e})_{(0)}=1$,
that is, a vertex operator algebra $V_{e}$ such that $V$ is a
subalgebra of $V_{e}$ and $\dim (V_{e})_{(0)}=1$.

\item A  haploid $\mathcal{C}$-algebra $V_{e}$ with
$\th_{V_{e}}=1_{V_{e}}$.
\end{enumerate}
\end{thm}
\pf Let $V_{e}$ be a vertex operator algebra such that $V$ is a
subalgebra of $V_{e}$.  Then $V_{e}$ is a $V$-module.  Thus it is also
an object of $\mathcal{C}$.  Since $V_{e}$ is a vertex operator
algebra, we have the vertex operator map $Y_{e}$ for $V_{e}$. This
vertex operator map can be viewed as an intertwining operator for $V$
of type ${V_{e}} \choose {V_{e}V_{e}}$. Let $\mu: V_{e}\boxtimes
V_{e}\to V_{e}$ be the module map corresponding to the intertwining
operator $Y_{e}$. Also, since $V$ is a subalgebra of $V_{e}$, we have
a morphism $\iota_{V_{e}}: V\to V_{e}$.  We now verify that the triple
$(V_{e}, \mu, \iota_{V_{e}})$ is indeed a haploid $\mathcal{C}$-algebra.

We first prove the associativity.  {}From the construction, we have
\[
\overline{\mu}(u\boxtimes v)=Y_{e}(u, 1)v
\]
for $u, v\in V_{e}$, where 
$\overline{\mu}: \overline{V_{e}\boxtimes V_{e}}\to \overline{V}_{e}$
is the natural extension of $\mu$. Let $\mu_{z}$ be the morphism from 
$V_{e}\boxtimes_{P(z)}V_{e}$ to $V_{e}$ corresponding to the intertwining 
operator $Y_{e}$. Then we have $\mu=\mu_{1}$ and 
\[
\overline{\mu}_{z}(u\boxtimes v)=Y_{e}(u, z)v.
\]
Thus for $u, v, w\in V_{e}$ and $z_{1}, z_{2}$ satisfying
$|z_{1}|>|z_{2}|>|z_{1}-z_{2}|>0$, 
\begin{eqnarray}
\overline{\mu_{z_{1}} \circ(1_{V_e} \boxtimes \mu_{z_{2}})}
(u\boxtimes_{P(z_{1})}
(v\boxtimes_{P(z_{2})} w))&=&Y_{e}(u, z_{1})Y_{e}(v, z_{2})w\label{prod}\\
\overline{\mu_{z_{1}} \circ(\mu_{z_{1}-z_{2}} \boxtimes 1_{V_e})}
((u\boxtimes_{P(z_{1}-z_{2})}
v)\boxtimes_{P(z_{2})} w)&=&Y_{e}(Y_{e}(u, z_{1}-z_{2})v, z_{2})w.\nn
\label{iter}
\end{eqnarray}
By the associativity for $Y_{e}$,
\begin{equation}
Y_{e}(u, z_{1})Y_{e}(v, z_{2})w
=Y_{e}(Y_{e}(u, z_{1}-z_{2})v, z_{2})w.\label{assoc}
\end{equation}
We also have the associativity isomorphism 
\[
\mathcal{A}^{P(z_{1}-z_{2}), P(z_{2})}_{P(z_{1}), P(z_{2})}: 
V_{e}\boxtimes_{P(z_{1})} (V_{e}\boxtimes_{P(z_{2})} V_{e})\to 
(V_{e}\boxtimes_{P(z_{1}-z_{2})} V_{e})\boxtimes_{P(z_{2})} V_{e},
\]
characterized by 
\begin{equation}
\overline{\mathcal{A}^{P(z_{1}-z_{2}), P(z_{2})}_{P(z_{1}), P(z_{2})}}
(u\boxtimes_{P(z_{1})} (v\boxtimes_{P(z_{2})} w))
=(u\boxtimes_{P(z_{1}-z_{2})} v)\boxtimes_{P(z_{2})} w\label{assoc-iso}
\end{equation}
for $u,v, w\in V_{e}$.
Combining (\ref{prod})--(\ref{assoc-iso}), we obtain
\begin{equation}
(\mu_{z_{1}} \circ(1_{V_e} \boxtimes_{P(z_{1})} \mu_{z_{2}}))=
(\mu_{z_{2}} \circ(\mu_{z_{1}-z_{2}} \boxtimes_{P(z_{2})} 1_{V_e}))\circ
\mathcal{A}^{P(z_{1}-z_{2}), P(z_{2})}_{P(z_{1}), P(z_{2})}.\label{mu-assoc1}
\end{equation}
{}From (\ref{mu-assoc1}), we obtain
\begin{eqnarray}
\lefteqn{(\mu_{z_{1}} \circ(1_{V_e} \boxtimes_{P(z_{1})} \mu_{z_{2}}))
\circ (1_{V_e} \boxtimes_{P(z_{1})} 
\mathcal{T}_{\gamma_{2}})\circ \mathcal{T}_{\gamma_{1}}}\nn
&&=
(\mu_{z_{2}} \circ(\mu_{z_{1}-z_{2}} \boxtimes_{P(z_{2})} 1_{V_e}))\circ
\mathcal{A}^{P(z_{1}-z_{2}), P(z_{2})}_{P(z_{1}), P(z_{2})}
\circ (1_{V_e} \boxtimes_{P(z_{1})} 
\mathcal{T}_{\gamma_{2}})\circ \mathcal{T}_{\gamma_{1}},\label{mu-assoc2}\nn
\end{eqnarray}
where $z_{1}, z_{2}$ are real numbers satisfying
$z_{1}>z_{2}>z_{1}-z_{2}>0$, $\gamma_{1}$ and $\gamma_{2}$ are paths
in the real line from $1$ to $z_{1}$ and $z_{2}$, respectively, and
$\mathcal{T}_{\gamma_{1}}$ and $\mathcal{T}_{\gamma_{2}}$ are the
parallel transport isomorphisms associated to $\gamma_{1}$ and
$\gamma_{2}$, respectively.

{}From the construction, we have
\begin{equation}
(\mu_{z_{1}} \circ(1_{V_e} \boxtimes_{P(z_{1})} \mu_{z_{2}}))\circ (1_{V_e}
\boxtimes_{P(z_{1})} \mathcal{T}_{\gamma_{2}})\circ
\mathcal{T}_{\gamma_{1}} =\mu \circ(1_{V_e} \boxtimes
\mu).\label{mu-assoc3}
\end{equation}
Similarly,
\begin{equation}
(\mu_{z_{2}} \circ(\mu_{z_{1}-z_{2}} \boxtimes_{P(z_{2})} 1_{V_e}))\circ
(\mathcal{T}_{\gamma_{3}}\circ (\mathcal{T}_{\gamma_{4}}
\boxtimes_{P(z_{2})} 1_{V_e}))^{-1}
=(\mu \circ(\mu \boxtimes 1_{V_e})),\label{mu-assoc4}
\end{equation}
where $\gamma_{3}$ and $\gamma_{4}$ are paths in the real line from
$z_{2}$ and $z_{1}-z_{2}$ to $1$, respectively, and
$\mathcal{T}_{\gamma_{3}}$ and $\mathcal{T}_{\gamma_{4}}$ are the
parallel transport isomorphisms associated to $\gamma_{3}$ and
$\gamma_{4}$, respectively.  Combining
(\ref{mu-assoc2})--(\ref{mu-assoc4}) and (\ref{mu-assoc5}), we obtain
the associativity
\[
\mu \circ(1_{V_e} \boxtimes \mu)=(\mu \circ(\mu \boxtimes 1_{V_{e}}))\circ 
\mathcal{A}.
\]

Next we prove the commutativity. We recall that the braiding
isomorphism $\mathcal{R}$ is characterized by
\begin{equation}
\overline{\mathcal{R}}(u\boxtimes v)=e^{L(-1)}
\overline{\mathcal{T}}_{\gamma_{1}^{-}}
(v\boxtimes_{P(-1)} u)\label{commu-iso}
\end{equation}
where $u, v\in V_{e}$, ${\gamma_{1}^{-}}$ is a (clockwise) path from $-1$ to $1$ in the
closed upper half plane with $0$ deleted, $\mathcal{T}_{\gamma_{1}^{-}}$
is the corresponding parallel transport isomorphism and
$\overline{\mathcal{T}}_{\gamma_{1}^{-}}$ is the natural extension of
$\mathcal{T}_{\gamma_{1}^{-}}$ to the formal completion
$\overline{V_{e}\boxtimes V_{e}}$ of $V_{e}\boxtimes V_{e}$.  Thus
\begin{equation}
\overline{\mu}(\overline{\mathcal{R}}(u\boxtimes v))=
\overline{\mu}(e^{L(-1)}\mathcal{T}_{\gamma_{1}^{-}}
(v\boxtimes_{P(-1)} u)). \label{mu-comm1}
\end{equation}
Let $\mathcal{Y}$ be the intertwining operator corresponding to 
the intertwining map $\boxtimes_{P(1)}$ of type 
${V_{e}\boxtimes V_{e}\choose V_{e}\;V_{e}}$. Then
\begin{equation}
e^{L(-1)}\mathcal{T}_{\gamma_{1}^{-}}
(v\boxtimes_{P(-1)} u)=e^{L(-1)}\mathcal{Y}(v, e^{\pi i})u.
\label{mu-comm2}
\end{equation}
Since $\mu$ is a morphism, $\mu \circ \mathcal{Y}=Y_{e}$, and 
$Y_{e}$ satisfies the skew-symmetry
\[
Y_{e}(u, 1)v=e^{L(-1)}Y_{e}(v, -1)u,
\]
we have
\begin{eqnarray}
\mu(e^{L(-1)}\mathcal{Y}(v, e^{\pi i})u)&=&e^{L(-1)}Y_{e}(v, -1)u\nn
&=&Y_{e}(u, 1)v\nn
&=&\overline{\mu}(u\boxtimes v).
\label{mu-comm3}
\end{eqnarray}
Combining (\ref{mu-comm1})--(\ref{mu-comm3}), we obtain
\[
\overline{\mu}(\overline{\mathcal{R}}(u\boxtimes v))
=\overline{\mu}(u\boxtimes v)
\]
for $u, v\in V_{e}$,  or equivalently,
the commutativity
\[
\mu\circ \mathcal{R}=\mu.
\]

For the unit property, we note that the left unit isomorphism
$l_{V_{e}}: V_{e}\to V\boxtimes V_{e}$ is defined by
$l_{V_{e}}(u)=\mathbf{1}\boxtimes u$ and thus
\begin{eqnarray*}
(\mu\circ (\iota_{V_{e}}\boxtimes 1_{V_{e}})\circ l_{V_{e}})(u)
&=&\mu((\iota_{V_{e}}\boxtimes 1_{V_{e}})(\mathbf{1}\boxtimes u))\nn
&=&\mu(\mathbf{1}\boxtimes u)\nn
&=&Y_{e}(\mathbf{1}, 1)u\nn
&=&1_{V_{e}}(u)
\end{eqnarray*}
for $u\in V_{e}$.

Finally, we prove the uniqueness of the unit.  Let $f\in
\hom_{\mathcal{C}}(V, V_{e})$. Since $f$ preserves the grading and
$\dim (V_{e})_{(0)}=1$, it is clear that $f$ maps $\mathbf{1}$ to a
scalar multiple of $\mathbf{1}$. Since as a module $V$ is generated by
$\mathbf{1}$, $f$ is determined by the scalar. Conversely, given any
scalar, we can construct an element of $\hom_{\mathcal{C}}(V, V_{e})$
such that it maps $\mathbf{1}$ to the scalar times $\mathbf{1}$.  Thus
$\dim \hom_{\mathcal{C}}(V, V_{e})=1$.

Conversely, let $(V_{e}, \mu, \iota_{V_{e}})$ be a
$\mathcal{C}$-algebra.  In particular, $V_{e}$ is a $V$-module. The
module map $\mu: V_{e}\boxtimes V_{e}\to V_{e}$ corresponds to an
intertwining operator $Y_{e}$ of type ${V_{e}\choose V_{e}V_{e}}$ such
that
\begin{equation}
\overline{\mu}(u\boxtimes v)=Y_{e}(u, 1)v\label{converse1}
\end{equation}
for $u, v\in V_{e}$.  Since we have an injective morphism $\iota: V\to
V_{e}$, we can view the vacuum vector $\mathbf{1}$ and the conformal
element $\omega$ of $V$ as elements of $V_{e}$. We now show that
$(V_{e}, Y_{e}, \mathbf{1}, \omega)$ is a vertex operator algebra
satisfying Conditions 1, 2, 3, 4 and 5.

Since $\th_{V_{e}}=1_{V_{e}}$, $V_{e}$ is $\mathbb{Z}$-graded.  One
immediate consequence is that $Y_{e}(u, x)v\in V_{e}((x))$ for $u,
v\in V_{e}$.  The skew symmetry for $Y_{e}$ now follows immediately
{}from $\mu\circ \mathcal{R}=\mu$ and the vacuum property follows
immediately from the unit property $\mu\circ (\iota_{V_{e}}\boxtimes
1_{V_{e}})\circ l_{V_{e}}=1_{V_{e}}$.  The creation property follows
{}from the vacuum property and the skew-symmetry. The Virasoro algebra
relations and the $L(0)$-grading property follow from the fact that
$V_{e}$ is a $V$-module. The $L(-1)$-derivative property follows from
the fact that $Y_{e}$ is an intertwining operator.

We now prove the associativity.  As above, for any nonzero complex
number $z$, let $\mu_{z}: V_{e}\boxtimes_{P(z)}V_{e} \to V_{e}$ be the
module map corresponding to the intertwining operator $Y_{e}$.  By
definition, we have
\begin{equation}
\mu_{z}(u\boxtimes_{P(z)} v)=Y_{e}(u, z)v=(\mu\circ \mathcal{T}_{\gamma})
(u\boxtimes_{P(z)} v)\label{converse3}
\end{equation}
for $u, v\in V_{e}$, where $z$ is a nonzero complex number and
$\gamma$ is a path from $z$ to $1$ in the complex plane with a cut
along the positive real line. Also,
(\ref{mu-assoc3})--(\ref{mu-assoc4}) hold.

Compose both sides of the associativity 
\[
\mu \circ(1_{V_{e}} \boxtimes \mu)=(\mu \circ(\mu \boxtimes 1_{V_{e}}))\circ
\mathcal{A}
\]
for the $\mathcal{C}$-algebra $V_{e}$ with
\[
((1_{V_{e}} \boxtimes_{P(z_{1})} 
\mathcal{T}_{\gamma_{2}})\circ \mathcal{T}_{\gamma_{1}})^{-1},
\]
where $\gamma_{1}$ and $\gamma_{2}$ are paths from $1$ to $z_{1}$ and
$z_{2}$, respectively, in the complex plane with a cut along the
positive real line and $z_{1}$ and $z_{2}$ are complex numbers
satisfying $|z_{1}|>|z_{2}|>|z_{1}-z_{2}|>0$.  Then we obtain
\begin{eqnarray}
\lefteqn{\mu \circ(1_{V_{e}} \boxtimes \mu)\circ
((1_{V_{e}} \boxtimes_{P(z_{1})} 
\mathcal{T}_{\gamma_{2}})\circ \mathcal{T}_{\gamma_{1}})^{-1}}\nn
&&=(\mu \circ(\mu \boxtimes 1_{V_{e}}))\circ 
\mathcal{A}\circ ((1_{V_{e}} \boxtimes_{P(z_{1})} 
\mathcal{T}_{\gamma_{2}})\circ \mathcal{T}_{\gamma_{1}})^{-1}.\label{converse4}
\end{eqnarray}
Using (\ref{mu-assoc3})--(\ref{mu-assoc5}) and (\ref{converse4}),
we have
\begin{eqnarray}
\lefteqn{\mu_{z_{1}} \circ(1_{V_{e}} \boxtimes_{P(z_{1})} \mu_{z_{2}})}\nn
&&=\mu \circ(1_{V_{e}} \boxtimes \mu)\circ
((1_{V_{e}} \boxtimes_{P(z_{1})} 
\mathcal{T}_{\gamma_{2}})\circ \mathcal{T}_{\gamma_{1}})^{-1}\nn
&&=(\mu \circ(\mu \boxtimes 1_{V_{e}}))\circ 
\mathcal{A}\circ ((1_{V_{e}} \boxtimes_{P(z_{1})} 
\mathcal{T}_{\gamma_{2}})\circ \mathcal{T}_{\gamma_{1}})^{-1}\nn
&&=(\mu_{z_{2}} \circ(\mu_{z_{1}-z_{2}} \boxtimes_{P(z_{2})} 1_{V_{e}}))\circ
(\mathcal{T}_{\gamma_{3}}\circ (\mathcal{T}_{\gamma_{4}}
\boxtimes_{P(z_{2})} 1_{V_{e}}))^{-1}\nn
&&\quad \quad\quad\quad\quad\circ
\mathcal{A}\circ ((1_{V_{e}} \boxtimes_{P(z_{1})} 
\mathcal{T}_{\gamma_{2}})\circ \mathcal{T}_{\gamma_{1}})^{-1}\nn
&&=(\mu_{z_{2}} \circ(\mu_{z_{1}-z_{2}} \boxtimes_{P(z_{2})} 1_{V_{e}}))
\circ \mathcal{A}^{P(z_{1}-z_{2}), P(z_{2})}_{P(z_{1}), P(z_{2})}.
\label{converse5}
\end{eqnarray}

For the next step, we need the convergence of products and iterates of
intertwining operators for $V$. Because of the convergence,
$\overline{\mu_{z_{1}}} \circ(1_{V_{e}} \boxtimes_{P(z_{1})}
\overline{\mu_{z_{2}}})$ is well defined and is equal to
$\overline{\mu_{z_{1}} \circ(1_{V_{e}} \boxtimes_{P(z_{1})}
\mu_{z_{2}})}$. Similarly, $\overline{\mu_{z_{1}}}
\circ(\overline{\mu_{z_{1}-z_{2}}} \boxtimes_{P(z_{2})} 1_{V_{e}})$ is
well defined and is equal to $\overline{\mu_{z_{1}}
\circ(\mu_{z_{1}-z_{2}} \boxtimes_{P(z_{2})} 1_{V_{e}})}$.  Thus
(\ref{converse5}) gives
\begin{eqnarray}
\lefteqn{\overline{\mu_{z_{1}}} \circ(1_{V_{e}} \boxtimes_{P(z_{1})} 
\overline{\mu_{z_{2}}})}\nn
&&=\overline{\mu_{z_{1}}} \circ(\overline{\mu_{z_{1}-z_{2}}}
\boxtimes_{P(z_{2})} 1_{V_{e}})\circ \overline{\mathcal{A}}^{P(z_{1}-z_{2}), 
P(z_{2})}_{P(z_{1}), P(z_{2})}.\label{converse6}
\end{eqnarray}
Applying both sides of (\ref{converse6}) to $u\boxtimes_{P(z_{1})}(v
\boxtimes_{P(z_{2})}w)$ for $u, v, w\in V_{e}$,
pairing the result with $v'\in V_{e}$ and using (\ref{converse3})
and 
\[
\overline{\mathcal{A}^{P(z_{1}-z_{2}), 
P(z_{2})}_{P(z_{1}), P(z_{2})}}(u\boxtimes_{P(z_{1})}(v
\boxtimes_{P(z_{2})}w))=(u\boxtimes_{P(z_{1}-z_{2})}v)
\boxtimes_{P(z_{2})}w,
\]
we obtain the associativity 
\[
\langle v', Y_{e}(u, z_{1})Y_{e}(v, z_{2})w\rangle
=\langle v', Y_{e}(Y_{e}(u, z_{1}-z_{2})v, z_{2})w\rangle
\]
for $u, v, w\in V_{e}$, $v'\in V_{e}'$ and $z_{1}, z_{2}\in \mathbb{C}$
satisfying $|z_{1}|>|z_{2}|>|z_{1}-z_{2}|>0$.

Since associativity and skew-symmetry imply commutativity, and
associativity and commutativity imply rationality (see \cite{H2} and
\cite{H3}), we have proved that $V_{e}$ is a vertex operator algebra.

{}From the uniqueness of unit, we have $\dim \hom_{\mathcal{C}}(V,
V_{e})=1$.  Assume that there is an element of $(V_{e})_{(0)}$ which
is not proportional to $\one$. Then this element generates a
$V$-submodule of the $V$-module $V_{e}$. Since every $V$-module is
completely reducible, this $V$-submodule is completely reducible. Thus
there exists an irreducible $V$-submodule of this $V$-submodule that
is generated by an element of $(V_{e})_{(0)}$ not proportional to
$\one$.  Since any irreducible $V$-module having a nonzero element of
weight $0$ must be equivalent to $V$, this irreducible $V$-submodule
is equivalent to $V$. But this $V$-submodule is not equal to
$\iota_{V_{e}}(V)\subset V_{e}$ since its generator of weight $0$ is
not proportional to $\one_{e}$. Thus $\dim \hom_{\mathcal{C}}(V,
V_{e})>1$, and this is a contradiction.  Hence $\dim (V_{e})_{(0)}=1$.  \epfv

\begin{rema}
{\rm From the proof, we see that the additional assumptions in Theorem
\ref{main} (that is, the assumptions that $V$ is simple, every $V$-module is
completely reducible, for any irreducible $V$-module $W$ the weights
of elements of $W$ are nonnegative real numbers and $W_{(0)}\ne 0$ if and only if
$W$ is equivalent to $V$ viewed as a $V$-module, and $\dim
V_{(0)}=1$) are needed only in the proof concerning the uniqueness of
unit. If we remove these additional assumptions, then the proof above
shows that an extension $V_{e}\supset V$ 
is equivalent to a (not necessarily haploid) $\mathcal{C}$-algebra. }
\end{rema}

\begin{thm}\label{main2}
Let $V, V_e$, and $\mathcal{C}$ be as in Theorem~\ref{main}.  Then the
category of modules for the vertex operator algebra $V_e$ is
isomorphic to the category ${\rm Rep}^0 V_e$ as defined in \cite{KO},
where $V_e$ is considered as an algebra in $\mathcal{C}$.
\end{thm}
\pf A $V_{e}$-module is clearly an object of ${\rm Rep}^0 V_e$.
Conversely, the proof of Theorem \ref{main} shows that the vertex
operator map for an object of ${\rm Rep}^0 V_e$ satisfies the
associativity and the skew-symmetry. Since for modules, the results
that associativity and skew-symmetry imply commutativity and that
associativity and commutativity imply rationality still hold (actually
more general results hold for intertwining operators; see \cite{H2}
and \cite{H3}), we see that this object is actually a $V_{e}$-module.
This correspondence between the two categories is an isomorphism of
categories. \epfv

\begin{thm}\label{main3}
Let $V, \mathcal{C}$ and $V_e$ be as in Theorem~\ref{main}. Assume in
addition that $V_e$ is simple as a $\mathcal{C}$-algebra, that there
are finitely many inequivalent irreducible $V$-modules, and that
$\mathcal{C}$ is rigid.  Then the fusion rules for $V$ and for $V_{e}$
are finite, every $V_{e}$-module is completely reducible, and $V_e$ as
a vertex operator algebra has finitely many inequivalent irreducible
modules.
\end{thm}
\pf Notice that under the assumptions, the category of $V$-modules has
a natural structure of rigid balanced braided tensor category. In
particular, we can use the results of Section 1 in \cite{KO}.

Since the tensor product of two $V$-modules as a $V$-module is a
direct sum of irreducible $V$-modules, no irreducible $V$-module can
have infinite multiplicity in the decomposition of this tensor
product.  Thus the fusion rules among irreducible $V$-modules are
finite and the fusion rules among general $V$-modules are also finite.
Since $V$ is a subalgebra of $V_{e}$, $V_{e}$-modules are $V$-modules.
Thus intertwining operators among $V_{e}$-modules are intertwining
operators among these $V_{e}$-modules as $V$-modules.  Since all the
fusion rules for $V$-modules are finite, all the fusion rules for
$V_{e}$ are finite.  By Theorem 3.3 in \cite{KO}, we also know that
any $V_{e}$-module is completely reducible.

We now prove that there are only finitely many inequivalent
irreducible $V_{e}$-modules.  Assume that there are infinitely many.
Since every $V$-module is completely reducible, any $V_{e}$-module as
a $V$-module is equivalent to a direct sum of irreducible
$V$-modules. Since there are only finitely many inequivalent
irreducible $V$-modules, every irreducible $V_{e}$-module as a
$V$-module is equivalent to a direct sum of copies of these
irreducible $V$-modules. Take any irreducible $V_{e}$-module $W$.
Then the multiplicity of each irreducible $V$-module in the
decomposition of $W$ as a $V$-module is finite.  Since there are
infinitely many irreducible $V_{e}$-modules, there must be an
irreducible $V_{e}$-module $W_{1}$ such that the multiplicity of at
least one irreducible $V$-module in the decomposition of $W_{1}$ as a
$V$-module is larger than the corresponding multiplicity in the
decomposition of $W$. Since $W$ and $W_{1}$ are inequivalent as
$V$-modules, they are also inequivalent as $V_{e}$-modules.  We can
continue this process to find infinitely many inequivalent irreducible
$V_{e}$-modules $W_{i}$ for $i\in \mathbb{Z}_{+}$ such that $W_{i}$ as
a $V$-module is equivalent to a $V$-submodule of $W_{i+1}$.  In
particular, $W$ as a $V$-module is equivalent to a $V$-submodule of
$W_{i}$ for $i\in \mathbb{Z}_{+}$.  For simplicity, we shall identify
$W$ with the equivalent $V$-submodules of $W_{i}$ for $i\in
\mathbb{Z}_{+}$.  Then all these infinitely many irreducible
$V_{e}$-modules are generated by $W$.

By the main result of \cite{P} (see Theorem 1.10 in \cite{KO}), the
category of $V_{e}$-module is a braided tensor category with the
tensor product bifunctor $\boxtimes^{V_{e}}$ defined as in the proof
of Theorem 1.5 in \cite{KO} (here we use the superscript $V_{e}$
rather than the subscript to avoid confusion with the notation
$\boxtimes_{P(z)}$ discussed in Section 2).  For $u\in V_{e}$ and
$w\in W$, $u\boxtimes w\in \overline{V_{e}\boxtimes W}$. By
definition, $V_{e}\boxtimes^{V_{e}} W$ is a quotient of
$V_{e}\boxtimes W$. In particular, $\overline{V_{e}\boxtimes^{V_{e}}
W}$ is a quotient of $\overline{V_{e}\boxtimes W}$.  We use
$u\boxtimes^{V_{e}} w$ to denote the coset of $u\boxtimes w$ in
$\overline{V_{e}\boxtimes^{V_{e}} W}$.  For $i\in \mathbb{Z}_{+}$,
there is a linear map from $V_{e}\boxtimes^{V_{e}} W$ to $W_{i}$ such
that its extension maps $u\boxtimes^{V_{e}} w$ to $Y_{W_{i}}(u,
1)w$. Clearly, this is a surjective $V_{e}$-module map from
$V_{e}\boxtimes^{V_{e}} W$ to $W_{i}$ for $i\in \mathbb{Z}_{+}$.  In
particular, $W_{i}$ for $i\in \mathbb{Z}_{+}$ are irreducible
quotients of $V_{e}\boxtimes^{V_{e}} W$ by $V_{e}$-submodules.  Since
every $V_{e}$-module is completely reducible, 
the $V_{e}$-module $V_{e}\boxtimes^{V_{e}} W$
is equivalent to a direct sum of irreducible $V_{e}$-modules. Since
$V$ has only finitely many irreducible modules and $V_{e}\boxtimes^{V_{e}} W$
is a $V_{e}$-module whose homogeneous subspaces are
finite-dimensional, $V_{e}\boxtimes^{V_{e}} W$ must be a finite direct sum of
irreducible $V_{e}$-modules. In fact, if it is an infinite direct sum,
then as a $V$-module it is also an infinite direct sum of irreducible
$V$-modules. Since there are only finitely many inequivalent
irreducible $V$-modules, any such infinite direct sum must have
infinite-dimensional homogeneous spaces, a contradiction.  Thus any
quotient $V_{e}$-module of $V_{e}\boxtimes^{V_{e}} W$ must also be equivalent
to a finite direct sum of the direct summands of $V_{e}\boxtimes^{V_{e}} W$.
In particular, each $W_{i}$ for $i\in \mathbb{Z}_{+}$ must be
equivalent to a direct summand of $V_{e}\boxtimes^{V_{e}} W$.  Since
there can only be finitely many such direct summands, we have a
contradiction, so that there must be only finitely many inequivalent
irreducible $V_{e}$-modules.  \epfv

Using Theorems \ref{main}--\ref{main3}, we obtain the following
result:

\begin{thm}\label{main4}
Let $V$ be a simple vertex operator algebra such whose module category
$\mathcal{C}$ has a natural vertex tensor category structure as
described in Section 2. We assume in addition the following: There are
finitely many inequivalent irreducible $V$-modules; every $V$-module
is completely reducible; for any irreducible $V$-module $W$, the
weights of elements of $W$ are nonnegative and $W_{(0)}\ne 0$ if and
only if $W$ is equivalent to $V$; $\dim V_{(0)}=1$; and the braided
tensor category structure on $\mathcal{C}$ is rigid.  Then the
following two notions are equivalent:

\begin{enumerate}
\item An extension $V_{e}\supset V$, where $V_e$ is a vertex operator
algebra such that $V$ is a subalgebra of $V_{e}$, its module category
has a natural braided tensor category structure as described in
Section 2, and it also satisfies the additional assumptions above for
$V$ with $V$ replaced by $V_{e}$, that is: There are finitely many
inequivalent irreducible $V_{e}$-modules; every $V_{e}$-module is
completely reducible; for any irreducible $V_{e}$-module $W$, the
weights of elements of $W$ are nonnegative and $W_{(0)}\ne 0$ if and
only if $W$ is equivalent to $V_{e}$; and $\dim (V_{e})_{(0)}=1$.

\item A haploid $\mathcal{C}$-algebra $V_{e}$ with
$\th_{V_{e}}=1_{V_{e}}$.
\end{enumerate}
In addition, the category of modules for such a vertex operator
algebra $V_e$ is isomorphic to the category ${\rm Rep}^0 V_e$ as
defined in \cite{KO}, where $V_e$ is the corresponding
$\mathcal{C}$-algebra.
\end{thm}
\pf By Theorem \ref{main}, a $\mathcal{C}$-algebra $V_{e}$ with
$\th_{V_{e}}=1_{V_{e}}$ is equivalent to a vertex operator algebra
$V_{e}$ such that $V$ is a subalgebra of $V_{e}$ and $\dim
(V_{e})_{(0)}=1$.  By Theorem \ref{main2}, the category of modules for
such a vertex operator algebra $V_e$ is isomorphic to the category
${\rm Rep}^0 V_e$.  By the main result of \cite{P} (see Theorem 1.10
in \cite{KO}), the category ${\rm Rep}^0 V_e$ has a natural braided
tensor category structure.

By Theorem \ref{main3}, every $V_{e}$-module is completely reducible
and $V_e$ as a vertex operator algebra has finitely many irreducible
modules.  Since any $V_{e}$-module is a $V$-module, for any
irreducible $V_{e}$-module $W$, the weights of elements of $W$ are
nonnegative.  If there is a $V_{e}$-module $W$ such that $W_{(0)}\ne
0$, then $W$ as a $V$-module must contain a $V$-submodule equivalent
to $V$. Since $W$ is irreducible as a $V_{e}$-module, $W$ is generated
as a $V_{e}$-module by this $V$-submodule.  Since $V_{e}$ is generated
as a $V_{e}$-module by $V$ and is also irreducible as a
$V_{e}$-module, $W$ is equivalent to $V_{e}$ as a $V_{e}$-module.
\epfv

\begin{rema}
{\rm Let $V$ and $\mathcal{C}$ be as in Theorem~\ref{main4}. Then by
Lemma 1.20 and Theorem 1.15 in \cite{KO}, a $\mathcal{C}$-algebra
$V_{e}$ with $\th_{V_{e}}=1_{V_{e}}$ is rigid if and only if $V_{e}$
as a $V_{e}$-module is irreducible. In this case, the braided tensor
category of $V_{e}$-modules when $V_{e}$ is viewed as a vertex
operator algebra is rigid.}
\end{rema}

\begin{rema}
{\rm In the theorems given in this section, we have to assume that the
category of $V$-modules has a natural structure of vertex tensor
category structure, as described in the preceding section.  The
existence of only a braided tensor category structure (even if it is
rigid and balanced) does not seem to be enough for these results to
hold. }
\end{rema}

\renewcommand{\theequation}{A.\arabic{equation}}
\renewcommand{\thethm}{A.\arabic{thm}}
\setcounter{equation}{0}
\setcounter{thm}{0}

\section*{Appendix}

Here we recall general theorems stating that under suitable
conditions, the module category for a vertex operator algebra has a
natural vertex tensor category structure and in particular braided
tensor category structure.
We also recall a general theorem stating that under suitable
conditions, this braided tensor category is rigid and is in fact a
modular tensor category.

We first state some conditions that we shall be considering here
(where we are using notation and terminology from
\cite{HL3}--\cite{HL7}, \cite{H1}):

\begin{enumerate}
\item Finite reductivity (called ``rationality'' in \cite{HL3}): Every
$V$-module is completely reducible, there exist only finitely many
inequivalent irreducible modules, and the spaces of intertwining
operators among triples of irreducible modules are all finite
dimensional. (Note that this finite reductivity is different from the
various notions of ``rationality'' in \cite{Z} and other works. Also
note that the first two conditions (every $V$-module is completely
reducible and there exist only finitely many inequivalent irreducible
modules) together with the finite dimensionality of the homogeneous
subspaces of $V$-modules imply that every $V$-module is of finite
length.)

\item The convergence and extension property for products: For any
$V$-modules $W_1$, $W_2$, $W_3$, $W_4$ and $W_{5}$ and any intertwining
operators ${\cal Y}_1$ and ${\cal Y}_2$ of types ${W_4}\choose {W_1W_{5}}$
and ${W_{5}}\choose {W_2W_3}$, respectively, there exists an integer $N$
(depending only on ${\cal Y}_1$ and ${\cal Y}_2$), and for any
homogeneous elements $w_{(1)}\in W_1$ and $w_{(2)}\in W_2$ and any
elements $w_{(3)}\in W_3$ and $w'_{(4)}\in W'_4$, there exist
$M\in{\mathbb N}$, $r_{k}, s_{k}\in {\mathbb R}$ ($k=1,\dots,M$), and
analytic functions $f_{k}(z)$ on $|z|<1$ ($k=1, \dots, M$) satisfying
\[
\wt w_{(1)}+\wt w_{(2)}+s_{k}>N, \;\;\;k=1, \dots, M,
\]
such that
\[
\langle w'_{(4)}, {\cal Y}_1(w_{(1)}, x_2) {\cal Y}_2(w_{(2)},
x_2)w_{(3)}\rangle \lbar_{x_1= z_1, \;x_2=z_2}
\]
is absolutely convergent when $|z_1|>|z_2|>0$ and can be analytically
extended to the multivalued analytic function
\[
\sum_{k=1}^{M}z_2^{r_{k}}(z_1-z_2)^{s_{k}}f_{k}\left(\frac{z_1-z_2}{z_2}\right)
\]
in the region $|z_2|>|z_1-z_2|>0$.

\item Convergence of products of intertwining operators: For any $n
\geq 3$, any $V$-modules $W_0, \dots, W_{n+1}$ and $\widetilde{W}_1,
\dots, \widetilde{W}_{n-1}$, any intertwining operators
\[
{\cal Y}_{1}, {\cal Y}_{2}, \dots, {\cal Y}_{n}
\]
of types
\[
{W_{0}\choose W_{1}\widetilde{W}_{1}}, {\widetilde{W}_{1}\choose
W_{2}\widetilde{W}_{2}}, \dots, {\widetilde{W}_{n-1}\choose
W_{n}W_{n+1}},
\]
respectively, and any $w_{(0)}'\in W_{0}'$, $w_{(1)}\in W_{1}, \dots,
w_{(n+1)}\in W_{n+1}$, the series
\begin{equation}
\langle w_{(0)}', {\cal Y}_{1}(w_{(1)}, x_{1})\cdots {\cal
Y}_{n}(w_{(n)}, x_{n})w_{(n+1)}\rangle|_{x_{j}=z_{j}, \; j=1, \dots, n}
\end{equation}
is absolutely convergent in the region $|z_{1}|>\cdots> |z_{n}|>0$ and
its sum can be analytically extended to a multivalued analytic
function on the region given by $z_{i}\ne 0$, $i=1, \dots, n$,
$z_{i}\ne z_{j}$, $i\ne j$, such that for any set of possible singular
points with either $z_{i}=0$, $z_{i}=\infty$ or $z_{i}= z_{j}$ for
$i\ne j$, this multivalued analytic function can be expanded near the
singularity as a series having the same form as the expansion near the
singular points of a solution of a system of differential equations
with regular singular points.

\item Every finitely-generated lower-bounded generalized $V$-module
is a module.

\item $V$ is of positive energy ($V_{(0)}=\mathbb{C}\one$ and
$V_{(n)}=0$ for $n<0$) and the $V$-module $V'$ contragredient to $V$
is equivalent to $V$.

\item Every grading-restricted 
generalized $V$-module is a (finite) direct sum of 
irreducible $V$-modules.

\item $V$ is $C_{2}$-cofinite, that is, $\dim V/C_{2}(V)<\infty$ where
$C_{2}(V)$ is the subspace of $V$ spanned by the elements of the form
$\res_{x}x^{-2}Y(u, x)v$ for $u, v\in V$.  (This condition was
introduced by Zhu in \cite{Z}, where it was called ``Condition $C$.'')
  
\end{enumerate}

The following theorem was proved in \cite{HL3}--\cite{HL7}, \cite{H1}:

\begin{thm}\label{vertextensorcat}
Let $V$ be a vertex operator algebra satisfying Conditions 1, 2, 3 and
4 above.  Then the category $\mathcal{C}$ of $V$-modules has a natural
structure of vertex tensor category and, in particular, $\mathcal{C}$
has a natural structure of braided tensor category.
\end{thm}

In \cite{H4}, the first author proved that if all $V$-modules satisfy
a condition called the $C_{1}$-cofiniteness condition, then the fusion
rules are all finite, and if in addition all the grading-restricted
generalized $V$-modules are completely reducible, then Conditions 2
and 3 above hold.  In particular, one can replace Conditions 2 and 3
in the theorem by these conditions. In \cite{Hu6}, the first author 
proved that if $V$ satisfies Conditions 5 and 7, then any finitely
generated lower-bounded generalized $V$-module is grading restricted. Thus
Conditions 5, 6 and 7 imply Condition 4. 
One consequence of these results is
that Conditions 5, 6 and 7 imply Conditions 1, 2, 3 and 4.  In
particular, for a vertex operator algebra $V$ satisfying Conditions 5,
6 and 7, the category of $V$-modules has a natural structure of
braided tensor category.

The following theorem was proved in \cite{H5}:

\begin{thm}
Let $V$ be a simple vertex operator algebra satisfying Conditions 5, 6 and 7
above.  Then the braided tensor category $\mathcal{C}$ of $V$-modules
is rigid and satisfies the nondegeneracy condition. In particular,
$\mathcal{C}$ is a modular tensor category.
\end{thm}

Theorem \ref{vertextensorcat} has been generalized in a number of
directions, including in particular non-semisimple module categories,
by Zhang and the first and third authors.  Since the proof of Theorem
\ref{main} requires only vertex tensor category structure on suitable
module categories, the conclusion of Theorem \ref{main} remains valid
in this greater generality.

\noindent {\small \sc Department of Mathematics, Rutgers University,
110 Frelinghuysen Rd., Piscataway, NJ 08854-8019}

\noindent {\em E-mail address}: {\tt yzhuang@math.rutgers.edu}

\vspace{1em}

\noindent {\small \sc Department of Mathematics, State University of 
New York at Stony Brook,
Stony Brook, NY 11794}

\noindent {\em E-mail address}: {\tt kirillov@math.sunysb.edu}

\vspace{1em}

\noindent {\small \sc Department of Mathematics, Rutgers University,
110 Frelinghuysen Rd., Piscataway, NJ 08854-8019}

\noindent {\em E-mail address}: {\tt lepowsky@math.rutgers.edu}

\end{document}